\documentclass[12pt]{amsart}
\usepackage{eepic}
\def\N{{\Bbb N}}

\newtheorem{Theorem}{Theorem}[section]

\newtheorem{Lemma}[Theorem]{Lemma}

\theoremstyle{definition}

\newtheorem{Example}{Example}
\newtheorem{Problem}{Problem}
\newtheorem*{Prob}{Problem}
\theoremstyle{remark}

\begin{document}
\sloppy
\title{A class of reflection rigid Coxeter systems}
\author{Tetsuya Hosaka} 
\address{Department of Mathematics, Utsunomiya University, 
Utsunomiya, 321-8505, Japan}
\date{October 16, 2004}
\email{hosaka@cc.utsunomiya-u.ac.jp}
\thanks{Partly supported by the Grant-in-Aid for Scientific Research, 
The Ministry of Education, Culture, Sports, Science and Technology, Japan, 
(No.~15740029).}
\maketitle
\begin{abstract}
In this paper, 
we give a class of reflection rigid Coxeter systems.
Let $(W,S)$ be a Coxeter system. 
Suppose that (1) for each $s,t\in S$ such that $m(s,t)$ is odd, 
$\{s,t\}$ is a maximal spherical subset of $S$, 
(2) there does not exist 
a three-points subset $\{s,t,u\}\subset S$ such that $m(s,t)$ and $m(t,u)$ are odd, 
and (3) for each $s,t\in S$ such that $m(s,t)$ is odd, 
the number of maximal spherical subsets of $S$ 
intersecting with $\{s,t\}$ is at most two, 
where $m(s,t)$ is the order of $st$ in the Coxeter group $W$.
Then we show that the Coxeter system $(W,S)$ is reflection rigid.
This is an extension of a result of 
N.~Brady, J.P.~McCammond, B.~M\"uhlherr and W.D.~Neumann.
\end{abstract}

\section{Introduction and preliminaries}

The purpose of this paper is to give a class of reflection rigid Coxeter systems.
A {\it Coxeter group} is a group $W$ having a presentation
$$\langle \,S \, | \, (st)^{m(s,t)}=1 \ \text{for}\ s,t \in S \,
\rangle,$$ 
where $S$ is a finite set and 
$m:S \times S \rightarrow \N \cup \{\infty\}$ is a function 
satisfying the following conditions:
\begin{enumerate}
\item[(1)] $m(s,t)=m(t,s)$ for any $s,t \in S$,
\item[(2)] $m(s,s)=1$ for any $s \in S$, and
\item[(3)] $m(s,t) \ge 2$ for any $s,t \in S$ such that $s\neq t$.
\end{enumerate}
The pair $(W,S)$ is called a {\it Coxeter system}.
For a Coxeter group $W$, a generating set $S'$ of $W$ is called 
a {\it Coxeter generating set for $W$} if $(W,S')$ is a Coxeter system.
In a Coxeter system $(W,S)$, 
the conjugates of elements of $S$ are called {\it reflections}. 
We note that 
the reflections depend on the Coxeter generating set $S$ and not just on 
the Coxeter group $W$. 
Let $(W,S)$ be a Coxeter system.
For a subset $T \subset S$, 
$W_T$ is defined as the subgroup of $W$ generated by $T$, 
and called a {\it parabolic subgroup}.
If $T$ is the empty set, then $W_T$ is the trivial group.
A subset $T\subset S$ is called a {\it spherical subset of $S$}, 
if the parabolic subgroup $W_T$ is finite.

A {\it diagram} is an undirected graph $\Gamma$ 
without loops or multiple edges 
with a map $\text{Edges}(\Gamma)\rightarrow\{2,3,4,\ldots\}$ 
which assigns an integer greater than $1$ to each of its edges. 
Since such diagrams are used to define Coxeter systems, 
they are called {\it Coxeter diagrams}. 

Let $(W,S)$ and $(W',S')$ be Coxeter systems. 
Two Coxeter systems $(W,S)$ and $(W',S')$ are 
said to be {\it isomorphic}, 
if there exists a bijection 
$\psi:S\rightarrow S'$ such that 
$$m(s,t)=m'(\psi(s),\psi(t))$$ 
for any $s,t \in S$, where 
$m(s,t)$ and $m'(s',t')$ are the orders of $st$ in $W$ 
and $s't'$ in $W'$, respectively.

In general, a Coxeter group does not always determine 
its Coxeter system up to isomorphism.
Indeed some counter-examples are known.

\begin{Example}[{\cite[p.38 Exercise~8]{Bo}}, \cite{BMMN}]\label{Ex:B}
It is known that for an odd number $k\ge 3$, 
the Coxeter groups defined 
by the diagrams in Figure~\ref{fig1} are 
isomorphic and $D_{2k}$.
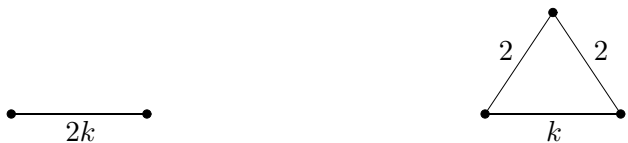
\begin{figure}[htbp]
\unitlength = 0.9mm
\begin{center}
\begin{picture}(80,28)(-40,-5)
\put(-25,0){\line(-1,0){20}}
\put(25,0){\line(1,0){20}}
\put(25,0){\line(2,3){10}}
\put(45,0){\line(-2,3){10}}
\put(-25,0){\circle*{1.3}}
\put(-45,0){\circle*{1.3}}
\put(25,0){\circle*{1.3}}
\put(45,0){\circle*{1.3}}
\put(35,15){\circle*{1.3}}
{\small
\put(-37,-4){$2k$}
\put(34,-4){$k$}
\put(27,8){$2$}
\put(41,8){$2$}
}
\end{picture}
\end{center}
\caption[Two distinct Coxeter diagrams for $D_{2k}$]{Two distinct Coxeter diagrams for $D_{2k}$}
\label{fig1}
\end{figure}
\end{Example}

\begin{Example}[\cite{BMMN}]\label{ex2}
It is known that
the Coxeter groups defined by 
the diagrams in Figure~\ref{fig2} are isomorphic by 
the {\it diagram twisting} (\cite[Definition~4.4]{BMMN}).
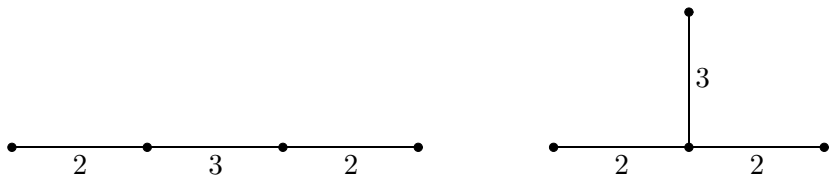
\begin{figure}[htbp]
\unitlength = 0.9mm
\begin{center}
\begin{picture}(140,30)(-70,-5)
\put(-5,0){\line(-1,0){60}}
\put(15,0){\line(1,0){40}}
\put(35,0){\line(0,1){20}}
\put(-5,0){\circle*{1.3}}
\put(-25,0){\circle*{1.3}}
\put(-45,0){\circle*{1.3}}
\put(-65,0){\circle*{1.3}}
\put(15,0){\circle*{1.3}}
\put(35,0){\circle*{1.3}}
\put(55,0){\circle*{1.3}}
\put(35,20){\circle*{1.3}}
{\small
\put(-16,-4){$2$}
\put(-36,-4){$3$}
\put(-56,-4){$2$}
\put(24,-4){$2$}
\put(44,-4){$2$}
\put(36,9){$3$}
}
\end{picture}
\end{center}
\caption[Coxeter diagrams for isomorphic Coxeter groups]{Coxeter diagrams for isomorphic Coxeter groups}
\label{fig2}
\end{figure}
\end{Example}

Here there exists the following natural problem.

\begin{Prob}[\cite{BMMN}, \cite{CD}]
When does a Coxeter group determine its Coxeter system up to isomorphism?
\end{Prob}

A Coxeter system $(W,S)$ is said to be {\it reflection rigid}, 
if for each Coxeter generating set $S'$ for $W$ such that 
the reflections of $(W,S)$ and $(W,S')$ coincide, 
$(W,S)$ and $(W,S')$ are isomorphic.

The following theorem was proved by 
N.~Brady, J.P.~McCammond, B.~M\"uhlherr and W.D.~Neumann in \cite{BMMN}.

\begin{Theorem}[{\cite[Theorem~3.9]{BMMN}}]
If $(W,S)$ is a Coxeter system such that 
$m(s,t)$ is even for each $s,t\in S$, 
then $(W,S)$ is reflection rigid, 
where $m(s,t)$ is the order of $st$ in the Coxeter group $W$.
\end{Theorem}

In this paper, we prove the following theorem 
which is an extension of the above theorem.

\begin{Theorem}\label{Thm}
Let $(W,S)$ be a Coxeter system.
Suppose that 
\begin{enumerate}
\item[(1)] for each $s,t\in S$ such that $m(s,t)$ is odd, 
$\{s,t\}$ is a maximal spherical subset of $S$, 
\item[(2)] there does not exist 
a three-points subset $\{s,t,u\}\subset S$ 
such that $m(s,t)$ and $m(t,u)$ are odd, and
\item[(3)] for each $s,t\in S$ such that $m(s,t)$ is odd, 
the number of maximal spherical subsets of $S$ 
intersecting with $\{s,t\}$ is at most two.
\end{enumerate}
Then $(W,S)$ is reflection rigid.
\end{Theorem}

Here we can not omit the condition (3) in this theorem by Example~\ref{ex2}.

\section{Proof of the theorem}

We first recall some basic results about Coxeter groups.

\begin{Lemma}[{\cite[p.12 Proposition~3]{Bo}}]\label{lem2-5}
Let $(W,S)$ be a Coxeter system and let $s,t\in S$.
Then $s$ is conjugate to $t$ if and only if 
there exists a sequence $s_1,\dots,s_n\in S$ such that 
$s_1=s$, $s_n=t$ and $m(s_i,s_{i+1})$ is odd for each $i\in\{1,\dots,n-1\}$.
\end{Lemma}

\begin{Lemma}[{\cite[Theorem~3.8]{BMMN}}]\label{lem3-1}
Let $(W,S)$ and $(W,S')$ be Coxeter systems.
If $R_S=R_{S'}$ then $|S|=|S'|$,
where $R_S$ and $R_{S'}$ are the sets 
of all reflections in $(W,S)$ and $(W,S')$, respectively.
\end{Lemma}

By Results~1.8, 1.9 and 1.10 in \cite{BMMN}, 
we obtain the following theorem.

\begin{Theorem}[cf.\ \cite{BMMN}]\label{thm4-0}
Let $(W,S)$ and $(W',S')$ be Coxeter systems. 
Suppose that there exists an isomorphism $\phi:W\rightarrow W'$. 
Then for each maximal spherical subset $T\subset S$, 
there exists a unique maximal spherical subset $T'\subset S'$ 
such that $\phi(W_T)=w'W'_{T'}(w')^{-1}$ for some $w'\in W'$.
\end{Theorem}

Concerning reflection rigidity of Coxeter systems, 
the following theorems have been proved.

\begin{Theorem}[{\cite[Theorem~3.9]{BMMN}}]\label{Thm:BMMN}
If $(W,S)$ is a Coxeter system such that 
$m(s,t)$ is even for any $s\neq t\in S$ 
(such $(W,S)$ is said to be {\it even}), 
then $(W,S)$ is reflection rigid.
\end{Theorem}

\begin{Theorem}[{\cite[Theorem~3.10]{BMMN}}]\label{Thm:finite}
If $(W,S)$ is a Coxeter system such that $W$ is finite, 
then $(W,S)$ is reflection rigid.
\end{Theorem}

Using these results, 
we prove the main theorem.

\begin{proof}[Proof of Theorem~\ref{Thm}]
Let $(W,S)$ be a Coxeter system.
Suppose that 
\begin{enumerate}
\item[(1)] for each $s,t\in S$ such that $m(s,t)$ is odd, 
$\{s,t\}$ is a maximal spherical subset of $S$, 
\item[(2)] there does not exist 
a three-points subset $\{s,t,u\}\subset S$ 
such that $m(s,t)$ and $m(t,u)$ are odd, and
\item[(3)] for each $s,t\in S$ such that $m(s,t)$ is odd, 
the number of maximal spherical subsets of $S$ 
intersecting with $\{s,t\}$ is at most two.
\end{enumerate}
Let $(W',S')$ be a Coxeter system.
We suppose that there exists an isomorphism 
$\phi:W\rightarrow W'$ such that 
for each reflection $r$ in $(W,S)$, 
$\phi(r)$ is a reflection in $(W',S')$.
Then we show that the Coxeter systems 
$(W,S)$ and $(W',S')$ are isomorphic.

We first show the following:
\begin{enumerate}
\item[$(1')$] for each $s',t'\in S'$ such that $m'(s',t')$ is odd, 
$\{s',t'\}$ is a maximal spherical subset of $S'$, 
\item[$(2')$] there does not exist 
a three-points subset $\{s',t',u'\}\subset S'$ 
such that $m'(s',t')$ and $m'(t',u')$ are odd, and
\item[$(3')$] for each $s',t'\in S'$ such that $m'(s',t')$ is odd, 
the number of maximal spherical subsets of $S'$ 
intersecting with $\{s',t'\}$ is at most two.
\end{enumerate}
Here $m'(s',t')$ is the order of $s't'$ in $W'$.

$(1')$ 
We immediately obtain $(1')$ from Theorems~\ref{thm4-0} and \ref{Thm:finite}.

$(2')$ 
Suppose that there exists 
a three-points subset $\{s',t',u'\}\subset S'$ 
such that $m'(s',t')$ and $m'(t',u')$ are odd.
Then $\{s',t'\}$ and $\{t',u'\}$ are maximal spherical subsets of $S'$ by $(1')$.
By Theorems~\ref{thm4-0}, 
there exist $s,t,u,v\in S$ 
such that $\phi^{-1}(W'_{\{s',t'\}})\sim W_{\{s,t\}}$ 
and $\phi^{-1}(W'_{\{t',u'\}})\sim W_{\{u,v\}}$, 
where $\{s,t\}\neq\{u,v\}$ (i.e., $|\{s,t,u,v\}|\ge 3$).
Then 
$$ s\sim t \sim \phi^{-1}(t')\sim\phi^{-1}(u')\sim u \sim v.$$
This contradicts (2) by Lemma~\ref{lem2-5}, since $|\{s,t,u,v\}|\ge 3$.
Thus $(2')$ holds.

$(3')$ 
Let $s',t'\in S'$ such that $m'(s',t')$ is odd.
Suppose that 
there exist maximal spherical subsets $T'_1$ and $T'_2$ of $S'$ 
intersecting with $\{s',t'\}$ such that 
the three sets $T'_1$, $T'_2$ and $\{s',t'\}$ are different.
We note that $\{s',t'\}$ is a maximal spherical subset of $S'$ by $(1')$.
By Theorems~\ref{thm4-0}, 
there exist $s,t\in S$ 
such that $\phi^{-1}(W'_{\{s',t'\}})\sim W_{\{s,t\}}$.
Also there exist maximal spherical subsets $T_1$ and $T_2$ 
such that $\phi^{-1}(W'_{T'_i})\sim W_{T_i}$ for each $i=1,2$.
Here $T_1$, $T_2$ and $\{s,t\}$ are different.
For each $i=1,2$, 
there exists an element $t_i\in T_i$ 
which is conjugate to $\phi^{-1}(T'_i\cap \{s',t'\})$.
Then $s,t,t_1,t_2$ are conjugate, 
and $\{t_1,t_2\}\subset \{s,t\}$ by (2).
This means that $t_i\in\{s,t\}\cap T_i\neq\emptyset$ for each $i=1,2$.
Hence $T_1$, $T_2$ and $\{s,t\}$ are 
different maximal spherical subsets of $S$ intersecting with $\{s,t\}$.
This contradicts (3).
Thus $(3')$ holds.

We note that $|S|=|S'|$ by Lemma~\ref{lem3-1}.

Let $s\in S$.
Since $s$ is a reflection in $(W,S)$, 
$\phi(s)$ is also a reflection in $(W',S')$ by the definition of $\phi$.
Then 
$$|\{a'\in S'\,|\, \phi(s)\sim a'\}|=1\ \text{or} \ 2$$
by Lemma~\ref{lem2-5}, $(1')$ and $(2')$.
Suppose that $|\{a'\in S'\,|\, \phi(s)\sim a'\}|=2$ and 
$\{a'\in S'\,|\, \phi(s)\sim a'\}=\{s',t'\}$.
Then $m'(s',t')$ is odd and $\{s',t'\}$ is a maximal spherical subset of $S'$.
If $T'$ is a maximal sherical subset of $S'$ 
such that $s'\in T'$ and $(W'_{T'},T')$ is even 
(such $T'$ is unique by $(3')$ if there exists), 
then $\{a'\in T'\,|\,\phi(s)\sim a'\}=\{s'\}$ by $(1')$ and $(2')$.

By the above argument, 
we can define a bijection $\psi:S\rightarrow S'$ as follows:
\begin{enumerate}
\item[(i)] $\phi(s)$ and $\psi(s)$ are conjugate for any $s\in S$.
\item[(ii)] Let $T$ be a maximal sherical subset of $S$ such that $(W_T,T)$ is even.
There exists a unique maximal sherical subset $T'$ of $S'$ 
such that $\phi(W_T)\sim W'_{T'}$ by Theorem~\ref{thm4-0}.
Then $\psi(s)\in T'$ for any $s\in T$.
\item[(iii)] Let $\{s,t\}$ be a maximal spherical subset of $S$
such that $m(s,t)$ is odd. 
There exists a unique maximal sherical subset $\{s',t'\}$ of $S'$ 
such that $\phi(W_{\{s,t\}})\sim W'_{\{s',t'\}}$ by Theorem~\ref{thm4-0}.
Then $\psi(\{s,t\})=\{s',t'\}$.
\end{enumerate}


We show that the bijection $\psi:S\rightarrow S'$ induces 
an isomorphism between $(W,S)$ and $(W',S')$.
Let $s,t\in S$ such that $m(s,t)$ is finite.
If $m(s,t)$ is odd, 
then $m(s,t)=m'(\psi(s),\psi(t))$ by (iii).
If $m(s,t)$ is even, 
then there exists a maximal spherical subset $T$ of $S$ 
such that $\{s,t\}\subset T$ and $(W_T,T)$ is even.
Then the restriction $\psi|T:T\rightarrow T'$ is a bijection such that 
$\phi(a)\sim\psi(a)$ for any $a\in T$.
By the proof of Theorem~\ref{Thm:BMMN} in \cite{BMMN}, 
$\psi|T:T\rightarrow T'$ induces an isomorphism between 
$(W_T,T)$ and $(W'_{T'},T')$.
Hence $m(s,t)=m'(\psi(s),\psi(t))$.
Thus $(W,S)$ and $(W',S')$ are isomorphic.
Therefore $(W,S)$ is reflection rigid.
\end{proof}

\section{Remarks}

The following theorem was proved by Bahls in \cite{B1}.

\begin{Theorem}[{\cite{B1}, cf.\ \cite[Theorem~5.2]{BM}}]\label{Thm:B1}
Let $(W,S)$ and $(W,S')$ be even Coxeter systems.
Then $(W,S)$ and $(W,S')$ are isomorphic.
\end{Theorem}

Here the following problem arises as an extension of Theorem~\ref{Thm:B1}.

\begin{Problem}
Let $(W,S)$ and $(W,S')$ be Coxeter systems such that 
\begin{enumerate}
\item[$(1)$] for any $s,t\in S$ such that $m(s,t)$ is odd, 
$\{s,t\}$ is a maximal spherical subset of $S$,
\item[(2)] there does not exist 
a three-points subset $\{s,t,u\}\subset S$ 
such that $m(s,t)$ and $m(t,u)$ are odd, 
\item[(3)] for each $s,t\in S$ such that $m(s,t)$ is odd, 
the number of maximal spherical subsets of $S$ 
intersecting with $\{s,t\}$ is at most two,
\item[$(1')$] for each $s',t'\in S'$ such that $m'(s',t')$ is odd, 
$\{s',t'\}$ is a maximal spherical subset of $S'$, 
\item[$(2')$] there does not exist 
a three-points subset $\{s',t',u'\}\subset S'$ 
such that $m'(s',t')$ and $m'(t',u')$ are odd, and
\item[$(3')$] for each $s',t'\in S'$ such that $m'(s',t')$ is odd, 
the number of maximal spherical subsets of $S'$ 
intersecting with $\{s',t'\}$ is at most two.
\end{enumerate}
Is it the case that $(W,S)$ and $(W,S')$ are isomorphic?
\end{Problem}

Recently, Bahls proved the following theorem in \cite{B2} 
which is an extension of a result in \cite{MW}.

\begin{Theorem}[\cite{B2}]\label{Thm:B2}
Let $(W,S)$ be a two-dimensional Coxeter system
(i.e.\ the Davis complex of $(W,S)$ is two-dimensional).
Then $(W,S)$ is reflection rigid up to diagram twisting.
\end{Theorem}

Here the following problem arises as an extension of this theorem.

\begin{Problem}
Let $(W,S)$ be a Coxeter system such that 
for each $s,t\in S$ if $m(s,t)$ is odd, 
then $\{s,t\}$ is a maximal spherical subset of $S$.
Is it the case that $(W,S)$ is reflection rigid up to diagram twisting?
\end{Problem}

The following theorem was proved in \cite{Ho}.

\begin{Theorem}[\cite{Ho}]\label{Thm:Ho}
Let $(W,S)$ and $(W,S')$ be two-dimensional Coxeter systems.
Then there exists $S''\subset W$ such that 
$(W,S'')$ is a Coxeter system which is isomorphic to $(W,S)$ 
and the sets of reflections in $(W,S'')$ and $(W,S')$ coincide.
\end{Theorem}

Theorems~\ref{Thm:B2} and \ref{Thm:Ho} implies the following theorem.

\begin{Theorem}\label{Thm:B2Ho}
Let $(W,S)$ and $(W,S')$ be two-dimensional Coxeter systems.
Then $(W,S)$ and $(W,S')$ are isomorphic up to diagram twisting.
\end{Theorem}

The following problem arises as an extension of this theorem.

\begin{Problem}
Let $(W,S)$ and $(W,S')$ be Coxeter systems such that 
\begin{enumerate}
\item[$(1)$] for any $s,t\in S$ such that $m(s,t)$ is odd, 
$\{s,t\}$ is a maximal spherical subset of $S$, and 
\item[$(1')$] for any $s',t'\in S'$ such that $m'(s',t')$ is odd, 
$\{s',t'\}$ is a maximal spherical subset of $S'$.
\end{enumerate}
Is it the case that $(W,S)$ and $(W,S')$ are isomorphic up to diagram twisting?
\end{Problem}

These problems are open.

%

%
\end{document}